\documentclass[reqno,a4paper]{amsart}
\usepackage{amssymb}%
\usepackage[hmarginratio=1:1]{geometry}%
\usepackage{ifpdf}
\usepackage{mathrsfs}
\usepackage{listings}
\ifpdf
 \usepackage[hyperindex]{hyperref}%,pagebackref
\else
 \expandafter\ifx\csname dvipdfm\endcsname\relax
 \usepackage[hypertex,hyperindex]{hyperref}
 \else
 \usepackage[dvipdfm,hyperindex]{hyperref}
 \fi
\fi
\usepackage{color}
\allowdisplaybreaks[4]
\theoremstyle{plain}
\newtheorem{theorem}{Theorem}[section]
\newtheorem{pro}{Proposition}[section]

\newtheorem{mr}{Monotonicity}[section]
\theoremstyle{remark}
\newtheorem{rem}{Remark}[section]
\theoremstyle{plain} 

\theoremstyle{definition} \newtheorem{definition}{Definition}[section]

\numberwithin{equation}{section}

\begin{document}
\title[]{$Y$-function and L'Hospital-type Monotonicity Rules with Nabla and Diamond-Alpha  Derivatives on Time Scales}
\author{Xiao-Yue Du, Zhong-Xuan Mao, Jing-Feng Tian*}

\address{Xiao-Yue Du \\
Hebei Key Laboratory of Physics and Energy Technology\\
Department of Mathematics and Physics\\
North China Electric Power University \\
Yonghua Street 619, 071003, Baoding, P. R. China}
\email{duxiaoyue\symbol{64}ncepu.edu.cn}

\address{Zhong-Xuan Mao\\
Hebei Key Laboratory of Physics and Energy Technology\\
Department of Mathematics and Physics\\
North China Electric Power University\\
Yonghua Street 619, 071003, Baoding, P. R. China}
\email{maozhongxuan000@gmail.com}

\address{Jing-Feng Tian\\
Hebei Key Laboratory of Physics and Energy Technology\\
Department of Mathematics and Physics\\
North China Electric Power University\\
Yonghua Street 619, 071003, Baoding, P. R. China}
\email{tianjf\symbol{64}ncepu.edu.cn}

\date{\today}
%criterion
\thanks{*Corresponding author: Jing-Feng Tian, e-mail addresses, tianjf\symbol{64}ncepu.edu.cn}
\maketitle
\begin{abstract}
The main objective of this paper is to establish the $Y$-function and L'Hospital-type monotonicity rules with nabla and diamond-alpha derivatives on time scales. \\
keywords: Monotonicity rule; L'Hospital-type monotonicity rule; Nabla derivatives; Diamond-Alpha derivatives; Time scale.
\end{abstract}

\section{Introduction}
If $-\infty \le a\le b\le + \infty $, if $f$ and $g$ are differentiable functions on $(a,b)$, and if $g^{\prime}\ne0$ on $(a,b)$, then Yang \cite{ref1} gave the following interesting auxiliary founction
\begin{equation}\label{1.1}
Y_{f,g}:=\frac{f ^{\prime}}{g ^{\prime}}g-f .
\end{equation}
 This function was explicitly named the $Y$-function by Tian et al. \cite{ref2}. Next, we present some properties of $Y$-function.
\begin{pro}\cite{ref1}
Let $-\infty \leq a<b\leq \infty $. Let $f$ and $g$ be
differentiable functions on $(a,b)$ and let $g^{\prime}\neq 0$ on $(a,b)$.
Let $Y_{f,g}$ be defined on $\left(a,b\right) $ by (\ref{1.1}). Then the following statements are true.\\
(i) The function is even with respect to $g$ and odd with respect to $f$, that is,
\begin{equation}
Y_{f,g}\left( x\right) =Y_{f,-g}\left( x\right) =-Y_{-f,g}\left( x\right)
=-Y_{-f,-g}\left( x\right).
\end{equation}
(ii) If $f$ and $g$ are twice differentiable on $(a,b)$, then
\begin{equation}
Y_{f,g}^{\prime}=\left( \frac{f^{\prime}}{g^{\prime}}\right) ^{\prime}g.
\end{equation}
(iii) If $g\neq 0$ on $\left( a,b\right) $, then
\begin{equation}
\left( \frac{f}{g}\right) ^{\prime}=\frac{g^{\prime}}{g^{2}}Y_{f,g}.
\end{equation}

\end{pro}

In 1982, Cheeger, Gromov and Taylor \cite{ref3} (see also \cite{M-T-2023})  presented the following  monotonicity rule for the ratio.

\begin{theorem}
If function $f$, $g$ are positive integrable on $\mathbb{R}$ and satisfy $f/g$ is decreasing, then the function
\begin{equation*}
x\mapsto \frac{\int_{0}^{x}f(t)dt}{\int_{0}^{x}g(t)dt}
\end{equation*}%
is decreasing.
\end{theorem}
In 1993, Anderson et al. \cite{ref4} (see also \cite{Qi-2022}) developed the following rule by using the Cauchy mean value theorem.
\begin{theorem} \label{Theorem1.2}
If differentiable functions $f,g$ defined on $(a,b]$ satisfy that $g^{\prime }(x)\neq 0$ for
all $x\in (a,b]$. If $x\mapsto f^{\prime }/g^{\prime }$ is increasing (decreasing)
on $(a,b]$, then
\begin{equation}
x\mapsto \frac{f(x)-f(a)}{g(x)-g(a)}
\end{equation}
is increasing (decreasing) on $(a,b]$.
\end{theorem}

Theorem \ref{Theorem1.2} has been proved by researchers. Different scholars have given it different names, we call this rule ``L'Hospital-type monotonicity rule'' in this paper \cite{ref5}.

\begin{theorem}
Let $f$ and $g$ be differentiable functions on the interval $(a,b)$, which satisfy the derivative $g^{\prime }\ne 0$ and does not change sign on $(a,b)$.
Suppose that $f(a^{+})=g(a^{+})=0$ or $f(b^{-})=g(b^{-})=0$. If $f^{\prime}/g^{\prime }$ is increasing (decreasing) on $(a,b)$ then $f/g$ is increasing (decreasing) on $(a,b)$.
\end{theorem}

Let $f(t)$ be diffrerntiable on time scales in the delta and nabla senses \cite{ref6}. For all $t\in\mathbb{T}$, Sheng et al. \cite{ref6} defined the diamond-alpha derivative and integral as follows
\begin{equation}
f^{\Diamond_{\alpha } }(t) =\alpha f^{\Delta } (t)+(1-\alpha )f^{\nabla}(t),\qquad  0\le \alpha \le 1
\end{equation}
and
\begin{equation*}
\int_{\beta }^{\gamma } f(t)\Diamond _{\alpha } t=\alpha\int_{\beta }^{\gamma } f(t)\Delta t+(1-\alpha )\int_{\beta }^{\gamma } f(t)\nabla t  ,\qquad  0\le \alpha \le 1.
\end{equation*}

We may notice that the diamond-alpha derivative reduces to the standard delta derivative as $\alpha=1$, or the standard nabla derivative as $\alpha=0$, while it respresents a weighted ``dynamic derivative" for $\alpha\in(0,1)$.

\section{${Y}$-function and L'Hospital-type Monotonicity Rules with Nabla Derivatives on Time Scales}
In this section, we will present ${Y}$-function and three monotonicity rules with nabla derivatives on time scales.

Before giving  Monotonicity rule \ref{Monotonicity rule 2.1}, we need to first present the so-called nabla $\mathcal{Y}$-function.
\begin{definition}
Let functions $\varphi , \psi$ defined on $[\alpha , \beta]_\mathbb{T}$ be nabla differential with $\psi^{\nabla}\not = 0$. Then the nabla $\mathcal{Y}$-function is defined by
\begin{equation}	
\mathcal{Y}_{\varphi,\psi }(x): = \frac{\varphi ^{\nabla }(x)}{\psi ^{\nabla}(x)}\psi (x)-\varphi(x),\quad x\in\left (\alpha ,\beta  \right ]_\mathbb{T}.
\end{equation}
\end{definition}
\begin{rem}
Clearly, if $\mathbb{T} = \mathbb{R}$, then nabla $\mathcal{Y}$-function reduces to $Y$-function
\begin{equation*}
Y_{\varphi ,\psi }(x): =
\frac{\varphi^{\prime}(x)}{\psi^{\prime}(x)}\psi(x) - \varphi(x) ,\quad x\in\left[\alpha , \beta\right].
\end{equation*}
\end{rem}
Based on the definition of the nabla $\mathcal{Y}$-function and simple calculation, we obtain the following interesting properties.

\begin{pro}\label{Propositions 2.1}
The following three assertions are valid.\\
(1) Nabla $\mathcal{Y}$-function has symmetry relations, namely,
\begin{equation}	
\mathcal{Y}_{\varphi , \psi} =
\mathcal{Y}_{\varphi , -\psi } =
\mathcal{Y}_{ -\varphi , \psi } =
-\mathcal{Y}_{ -\varphi , -\psi }.
\end{equation}
(2) If the function $\varphi^{\nabla}(s) / \psi^{\nabla}(s)$ is nabla differentiable, then the function $\mathcal{Y}_{\varphi,\psi }$ is nabla differential with
\begin{equation}\label{formula(2.3)}
\mathcal{Y}^{\nabla } _{\varphi, \psi} = (\frac{\varphi^{\nabla }}{\psi^{\nabla }})^{\nabla}\psi ^{\rho }.
\end{equation}
(3) If $\psi\not = 0$, then
\begin{equation}\label{formula(2.4)}
(\frac{\varphi }{\psi })^{\nabla } =
\frac{\psi ^{\nabla}}{\psi \psi ^{\rho }}\mathcal{Y}_{\varphi , \psi }.
\end{equation}
\end{pro}
\begin{pro}\label{Propositions 2.2}
Suppose that the function $\varphi^{\nabla} / \psi^{\nabla}$ is nabla differentiable, then the following two assertions are valid.
	
(1) The function $\mathcal{Y}_{\varphi, \psi}$ is increasing (decreasing) if $\psi^{\rho}\ge0$ and the function $\varphi^{\nabla}(s) / \psi^{\nabla}(s)$ is increasing (decreasing).
	
(2) The function $\mathcal{Y}_{\varphi, \psi }$ is decreasing (increasing) if $\psi^{\rho}\le0$ and the function $\varphi^{\nabla}(s) / \psi^{\nabla}(s)$ is increasing (decreasing).
\end{pro}
\begin{rem}
We have omitted here the proof of Propositions \ref{Propositions 2.1} and \ref{Propositions 2.2}.
\end{rem}
Then, we give the first kind of L'Hospital-type monotonicity rules with nabla derivatives on time scales.

\begin{mr}\label{Monotonicity rule 2.1}
Suppose that the functions $\varphi,\psi$ defined on $[\alpha,\beta]_{\mathbb{T}}$ satisfy that $\varphi^{\nabla}$ and $\psi^{\nabla}$ are rd-continuous as well as the sign of $\psi$ does not change. If the function $s\mapsto \varphi^{\nabla}(s) / \psi^{\nabla}(s)$ is increasing (decreasing) on $[\alpha,\beta]_\mathbb{T}$, then the functions
\begin{equation*}
s \mapsto \frac{\varphi(s) - \varphi(\alpha)}{\psi (s) - \psi(\alpha)} \qquad \ and \qquad s \mapsto \frac{\varphi(s) - \varphi(\beta)}{\psi (s) - \psi(\beta)}
\end{equation*}
are increasing (decreasing) on $[\alpha,\beta]_\mathbb{T}$.
\end{mr}

\begin{proof}
Since function $\varphi^{\nabla}$ and $\psi^{\nabla}$ are rd-continuous, we have	
\begin{equation*}
{\varphi(s) - \varphi (\alpha)} =
\int_{\alpha}^{s} \varphi^{\nabla}(\tau)\nabla\tau
\end{equation*}	
and
\begin{equation*}	
{\psi(s) - \psi (\alpha)} =
\int_{\alpha}^{s} \psi^{\nabla}(\tau) \nabla\tau.
\end{equation*}
By using formulas
\begin{equation*}
{\left( \frac{f(s)}{g(s)} \right)^{\nabla } =
\frac{f^{\nabla}(s) g(s) - f(s) g^{\nabla}(s)} {g(\rho(s))g(s)}}
\end{equation*}
and
\begin{equation*}
\left(\int_{\alpha }^{s} f(\tau)\nabla\tau \right)^{\nabla } = f(s) ,\quad
f\in C_{rd}(\mathbb{T}, \mathbb{R}),
\end{equation*}
we obtain
\begin{equation*}
\begin{aligned}
\nonumber&\int_{\alpha}^{s} \psi ^{\nabla }(\tau )\nabla\tau
\int_{\alpha}^{\rho(s) } \psi ^{\nabla}(\tau)\nabla\tau
\left(\frac{\varphi (s)-\varphi (\alpha )}{\psi (s)-\psi (\alpha )}\right)^{\nabla }\\
\nonumber& = \int_{\alpha}^{s} \psi ^{\nabla }(\tau )\nabla \tau \int_{\alpha}^{\rho(s) }\psi ^{\nabla}(\tau)\nabla\tau
\left(\frac{\int_{\alpha}^{s} \varphi^{\nabla }(\tau )\nabla\tau}{\int_{\alpha}^{s } \psi ^{\nabla }(\tau )\nabla\tau}\right)^{\nabla }\\
\nonumber& = \varphi ^{\nabla }(s) \int_{\alpha}^{s } \psi ^{\nabla }(\tau )\nabla\tau - \psi ^{\nabla }(s) \int_{\alpha}^{s }\varphi^{\nabla}(\tau)\nabla\tau\\
\nonumber& = \int_{\alpha }^{s}\psi^{\nabla }(\tau ) \psi ^{\nabla }(s)\left(\frac{\varphi ^{\nabla }(s)}{\psi ^{\nabla }(s)} - \frac{\varphi ^{\nabla }(\tau )}{\psi ^{\nabla }(\tau )}\right)\nabla\tau
\end{aligned}
\end{equation*}
which implies the desired result.
\end{proof}

\begin{rem}
The monotonicity rule \ref{Monotonicity rule 2.1} was given firstly by Martins and Torres in \cite{Martins}, and here we give
a simple proof of this rule.
\end{rem}

Based on the monotonicity of the function $\varphi^{\nabla}(s) / \psi^{\nabla}(s)$ and the value of $\mathcal{Y}_{\varphi , \psi }(b)$, where $b\in \left \{ \alpha ,\beta \right\}$, we can identify the monotonicity of the function $\varphi/\psi$. Now we present the second kind of L'Hospital-type monotonicity rules with nabla derivatives on time scales.

\begin{mr}
Let functions $\varphi$ and $\psi$ defined on $[\alpha , \beta]_\mathbb{T}$ be nabla differential with $\psi>0$ on $[\alpha , \beta]_\mathbb{T}$ and $\psi^{\nabla}>(<)0$ on $[\alpha , \beta)_\mathbb{T}$.

(1) If the function $\varphi^{\nabla} / \psi^{\nabla}$ is increasing (decreasing) and $\mathcal{Y}_{\varphi , \psi } (\alpha)(\mathcal{Y}_{\varphi , \psi}(\beta ))\ge 0$ then the function $\varphi/\psi$ is increasing (decreasing).

(2) If the function $\varphi^{\nabla} / \psi^{\nabla}$ is increasing (decreasing) and $\mathcal{Y}_{\varphi,\psi} (\beta)(y_{\varphi,\psi}(\alpha))\le0$, then the function  $\varphi/\psi$ is decreasing (increasing).

(3) If the function $\varphi^{\nabla}/\psi^{\nabla}$ is increasing and $\mathcal{Y}_{\varphi ,\psi}(\alpha)\le 0 ,\mathcal{Y}_{\varphi ,\psi }(\beta)\ge0$, then there exists a number $x_{0}\in{[\alpha,\beta]}$ such that the function $\varphi / \psi$ is decreasing (increasing) on $[\alpha,x_{0}]_\mathbb{T}$ and increasing (decreasing) on $[x_{0} , \beta]_\mathbb{T}$.

(4) If the function $\varphi^{\nabla} / \psi^{\nabla}$ is decreasing and $\mathcal{Y}_{\varphi,\psi}(\alpha)\ge 0,\mathcal{Y}_{\varphi ,\psi} (\beta)\le 0$, then there exists a number $x_{0}\in{[\alpha,\beta]}$ such that the function $\varphi / \psi$ is increasing (decreasing) on $[\alpha , x_{0}]_\mathbb{T}$ and decreasing (increasing) on $[x_{0},\beta]_\mathbb{T}$.
\end{mr}

\begin{proof}
Based on formula \ref{formula(2.4)}, we know that if both $\psi^{\nabla}$ and $\mathcal{Y}_{\varphi , \psi }$ are nonnegative or nonpositive, then the function $\varphi/\psi$ is increasing, if one of $\psi^{\nabla}$ and $\mathcal{Y}_{\varphi , \psi }$ is nonnegative and the other is nonpositive, then the function $\varphi/\psi$ is decreasing.

(1) If $\psi>0$ and the function $\varphi^{\nabla} / \psi^{\nabla}$ is increasing (decreasing), then the function $\mathcal{Y}_{\varphi , \psi}$ is increasing (decreasing) based on Proposition \ref{Propositions 2.2}. Together with the fact that $\mathcal{Y}_{\varphi , \psi}(\alpha )(\mathcal{Y}_{\varphi , \psi}(\beta))\ge0$, we obtain that the function $\mathcal{Y}_{\varphi ,\psi}$ is nonnegative. Thus we receive that the function $\varphi/\psi$ is increasing (decreasing).

(2) In the same method, we obtain that the function $\mathcal{Y}_{\varphi ,\psi}$ is increasing  (decreasing) and $\mathcal{Y}_{\varphi ,\psi }$ is non-positive. Thus the function $\varphi/\psi$ is decreasing (increasing).

(3) The function $\mathcal{Y}_{\varphi ,\psi}$ is increasing with with  $\mathcal{Y}_{\varphi,\psi}(\alpha )\le 0$ and $\mathcal{Y}_{\varphi ,\psi }(\beta )\ge 0$ deduces that there exists $x_{0}\in{[\alpha,\beta]}$ such that $\mathcal{Y}_{\varphi ,\psi}(x)\le 0$ for all $x\in \left [\alpha ,x_{0} \right] _\mathbb{T}$ and $\mathcal{Y}_{\varphi , \psi}(x)\le 0$ for all $x\in \left [x_{0} ,\beta \right] _\mathbb{T}$. Thus the function $\varphi/\psi$ decreasing (increasing) on $[\alpha , x_{0}]_\mathbb{T}$ and it is increasing (decreasing) on $[x_{0},\beta]_\mathbb{T}$, which is the desired conclusion.

(4) The function $\mathcal{Y}_{\varphi ,\psi}$ is decreasing with with $\mathcal{Y}_{\varphi,\psi}(\alpha )\ge 0$ and $\mathcal{Y}_{\varphi ,\psi} (\beta)\le 0$ deduces that there exists $x_{0}\in{[\alpha,\beta]}$ such that $\mathcal{Y}_{\varphi ,\psi}(x)\ge0$ for all $x\in \left [\alpha ,x_{0} \right]_\mathbb{T}$ and $\mathcal{Y}_{\varphi ,\psi}(x )\le0$ for all $x\in \left[ x_{0},\beta \right ] _\mathbb{T}$. Then we obtain the function $\varphi/\psi$ is increasing(decreasing) on $[\alpha,x_{0}]_\mathbb{T}$ and $\varphi/\psi$ is decreasing (increasing) on $[x_{0},\beta]_\mathbb{T}$, which is the desired conclusion.

Finally, we provided the third kind of L'Hospital-type monotonicity rules with nabla derivatives on time scales.
\end{proof}
\begin{mr}
Let functions $\varphi$ and $\psi$ defined on $[\alpha,\beta]_\mathbb{T}$ be nabla differential with $\varphi(\alpha )=\psi(\alpha)= 0$. If there exists a number $p\in \left (\alpha,\beta \right )$ such that the function $\varphi^{\nabla}/\psi^{\nabla}$ is increasing on $\left [\alpha ,p \right]_\mathbb{T}$ and it is increasing on $\left[p,\beta  \right]_\mathbb{T}$, then statements (i) and (ii) hold; if there exists a number $p\in \left (\alpha, \beta \right)$ such that the function $\varphi^{\nabla} / \psi^{\nabla}$ is decreasing on $\left[\alpha , p\right ]_\mathbb{T}$ and it is increasing on $\left[p,\beta  \right] _\mathbb{T}$, then statements (iii) and (iv) hold.

(i) If both $\frac{\psi ^{\nabla }}{\psi \psi ^{\rho}}$ and $\mathcal{Y}_{\varphi ,\psi}(\beta)$ are nonnegative or nonpositive, then the function $\varphi / \psi$ is increasing on $\left[ \alpha ,\beta \right]_\mathbb{T}$.

(ii) If one of $\frac{\psi^{\nabla }}{\psi \psi^{\rho }}$ and $\mathcal{Y}_{\varphi ,\psi }(\beta )$ is nonnegative and the other is nonpositive, then there exists a number $x_{p}$ such that the function $\varphi / \psi$ is increasing on $\left[\alpha,x_{p} \right]_\mathbb{T}$ and it is decreasing on $\left[ x_{p},\beta \right ]_\mathbb{T}$.

(iii) If both $\frac{\psi^{\nabla}}{\psi \psi ^{\rho}}$ and $\mathcal{Y}_{\varphi ,\psi}(\beta )$ are nonnegative or nonpositive, then there exists a number $x_{p}$ such that the function $\varphi / \psi$ is decreasing on $\left [\alpha ,x_{p} \right] _\mathbb{T}$ and it is increasing on $\left [x_{p},\beta \right]_\mathbb{T}$.

(iv) If one of $\frac{\psi^{\nabla }}{\psi \psi^{\rho }}$ and $\mathcal{Y}_{\varphi ,\psi }(\beta )$ is nonnegative and the other is nonpositive, then the function $\varphi/\psi$ is decreasing on $\left[\alpha,\beta \right]_\mathbb{T}$.
\end{mr}
\begin{proof}
Without prejudice to generality, we consider the situation that the function $\varphi^{\nabla} / \psi^{\nabla}$ is increasing and then decreasing. According to the formula \ref{formula(2.3)}, we obtain that the function $\mathcal{Y}_{\varphi ,\psi}$ is  increasing (decreasing) then decreasing (increasing) if $\psi ^{\rho } > (<) 0$. Then $\mathcal{Y}_{\varphi ,\psi }(\beta)$ is always nonnegative (nonpositive) if $\mathcal{Y}_{\varphi ,\psi }(\beta )\ge (\le )0$ and $\mathcal{Y}_{\varphi ,\psi}$ is first nonnegative (nonpositive) and then nonpositive (nonnegative) if $\mathcal{Y}_{\varphi,\psi }(\beta)\le (\ge)0$. Thus, the proof is completed based on the formula \ref{formula(2.4)}.
\end{proof}

\section{$Y$-function and L'Hospital-type Monotonicity Rules with Diamond-Alpha Derivatives on Time Scales}
In this section, the main purpose is to establish $Y$-function and the L'Hospital-type monotonicity rules with diamond-alpha derivatives on time scales.

Before giving  Monotonicity rule \ref{Monotonicity rule 3.1}, we present the so-called diamond-alpha $\mathcal{Y}$-function firstly.

\begin{definition}
Let functions $\varphi , \psi$ defined on $[\beta,\gamma]_\mathbb{T}$ be diamond-alpha differential with $\psi^{\Diamond_{\alpha}}\not = 0$. Then the diamond-alpha $\mathcal{Y}$-function is defined by
	
\begin{equation}	
\mathcal{Y}_{\varphi,\psi }(x): = \frac{\varphi ^{\Diamond_{\alpha}}(x)}{\psi ^{\Diamond_{\alpha}}(x)}\psi (x)-\varphi(x),\quad x\in\left [\beta,\gamma  \right ]_\mathbb{T}.
\end{equation}
	
\end{definition}

\begin{rem}
Clearly, if $\mathbb{T} = \mathbb{R}$, then diamond-alpha $\mathcal{Y}$-function reduces to $Y$-function
\begin{equation*}
Y_{\varphi ,\psi }(x): =
\frac{\varphi^{\prime}(x)}{\psi^{\prime}(x)}\psi(x) - \varphi(x) ,
\quad x\in\left[\beta,\gamma\right].
\end{equation*}
\end{rem}
Based on the definition of the diamond-alpha $\mathcal{Y}$-function and simple calculation, we obtain the following interesting properties.
\begin{pro}\label{Propositions 3.1}
The following three assertions are valid.\\
(i) Diamond-alpha $\mathcal{Y}$-function has symmetry relations, namely,
\begin{equation}	
\mathcal{Y}_{\varphi , \psi} =
\mathcal{Y}_{\varphi , -\psi } =
\mathcal{Y}_{ -\varphi , \psi } =
-\mathcal{Y}_{ -\varphi , -\psi }.
\end{equation}
(ii) If the function $\varphi^{\Diamond_{\alpha}}(x) / \psi^{\Diamond_{\alpha}}(x)$ is diamond-alpha differentiable and the function $\varphi(x),\psi(x)$ is delta and nabla differentiable, then the function $\mathcal{Y}_{\varphi,\psi }$ is diamond-alpha differential with
\begin{equation}\label{formula(1.3)}
\begin{aligned}
\mathcal{Y}^{\Diamond_{\alpha}} _{\varphi,\psi}
\nonumber&=\frac{\alpha}{\psi ^{\Diamond _{\alpha } }(x)\psi ^{\sigma \Diamond _{\alpha } }(x)}\Big(\big(\alpha \varphi ^{\Delta\Delta }(x)+(1-\alpha )\varphi ^{\nabla \Delta }(x))\psi ^{\sigma }(x)+\varphi ^{\Diamond _{\alpha }}(x) \psi ^{\Delta }(x)\big)\psi ^{\Diamond _{\alpha } }(x)\\
\nonumber&\quad-\big(\alpha\psi ^{\Delta \Delta }(x)+(1-\alpha )\psi ^{\nabla \Delta }(x)\big)\varphi ^{\Diamond _{\alpha } }(x) \psi(x)-\varphi ^{\Delta }(x)\psi ^{\Diamond _{\alpha } }(x)\psi ^{\sigma \Diamond _{\alpha } }(x)\Big)\\
\nonumber&\quad+\frac{1-\alpha}{{\psi ^{\Diamond _{\alpha } }(x)\psi ^{\rho\Diamond _{\alpha } }(x)}} \Big(\big(\alpha \varphi ^{\Delta\nabla }(x)+(1-\alpha )\varphi ^{\nabla \nabla }(x))\psi ^{\rho  }(x)+\psi ^{\Diamond _{\alpha } }(x)\psi ^{\nabla   }(x)\big)\psi ^{\Diamond _{\alpha } }(x)\\
\nonumber&\quad-\big(\alpha\psi ^{\Delta \nabla }(x)+(1-\alpha )\psi ^{\nabla \nabla }(x)\big)\varphi  ^{\Diamond _{\alpha } }(x)\psi(x)-\varphi ^{\nabla }(x)\psi ^{\Diamond _{\alpha } }(x)\psi ^{\rho\Diamond _{\alpha } }(x) \Big),\\
\nonumber&\qquad0\le \alpha\le1.
\end{aligned}
\end{equation}
(iii) If $\psi (x)\psi ^{\sigma } (x)\psi ^{\rho } (x)\ne 0$, and $\varphi /\psi :\mathbb{T}\to R$ is diamond-alpha differentiable at $x\in\mathbb{T}$, then
\begin{equation*}\label{(3.1)}
(\frac{\varphi }{\psi } )^{\Diamond_{\alpha}}(x)=\frac{(\varphi ^{\Diamond _{\alpha } } (x)\psi ^{\sigma } (x)\psi ^{\rho } (x)-\alpha \varphi ^{\sigma }(x)\psi ^{\rho }(x)\psi ^{\Delta }(x)
-(1-\alpha )\varphi ^{\rho }(x)\psi ^{\sigma }(x)\psi ^{\nabla }(x))}{\psi (x)\psi ^{\sigma } (x)\psi ^{\rho } (x)}
\end{equation*}
\end{pro}
\begin{proof}
We only need to show (ii) since proofs of (i) is straightforward from the definition of the diamond-alpha and (iii) has been proofed by \cite{ref6}. (ii) is obtained by using following formula
	
\begin{equation}
\begin{aligned}
\mathcal{Y}^{\Diamond} _{\varphi,\psi}
\nonumber&=\big(\frac{\varphi ^{\Diamond_{\alpha}}(x)}{\psi ^{\Diamond_{\alpha}}(x)}\psi(x)-\varphi(x)\big)^{\Diamond _{\alpha}}\\
\nonumber&=\alpha\big(\frac{\varphi ^{\Diamond_{\alpha}}(x)}{\psi ^{\Diamond_{\alpha}}(x)}\psi(x)-\varphi(x)\big)^{\Delta }
+(1-\alpha )\big(\frac{\varphi ^{\Diamond_{\alpha}}(x)}{\psi ^{\Diamond_{\alpha}}(x)}\psi(x)-\varphi(x)\big)^{\nabla }\\
\nonumber&=\alpha\big(\frac{\alpha \varphi ^{\Delta}(x)+
(1-\alpha)\varphi ^{\nabla }(x)}{\alpha\psi ^{\Delta}(x)+
(1-\alpha)\psi^{\nabla }(x)}\psi(x)-\varphi(x)\big)^{\Delta }\\
&\quad+(1-\alpha)\big(\frac{\alpha \varphi ^{\Delta}(x)+(1-\alpha)\varphi ^{\nabla }(x)}{\alpha\psi  ^{\Delta}(x)+(1-\alpha)\psi^{\nabla }(x)}\psi(x)-\varphi(x)\big)^{\nabla }\\
\nonumber&=\frac{\alpha}{\psi ^{\Diamond _{\alpha } }(x)\psi ^{\sigma \Diamond _{\alpha } }(x)}\Big(\big(\alpha \varphi ^{\Delta\Delta }(x)+(1-\alpha )\varphi ^{\nabla \Delta }(x))\psi ^{\sigma }(x)+\varphi ^{\Diamond _{\alpha }}(x) \psi ^{\Delta }(x)\big)\psi ^{\Diamond _{\alpha } }(x)\\
\nonumber&\quad-\big(\alpha\psi ^{\Delta \Delta }(x)+(1-\alpha )\psi ^{\nabla \Delta }(x)\big)\varphi ^{\Diamond _{\alpha } }(x) \psi(x)-\varphi ^{\Delta }(x)\psi ^{\Diamond _{\alpha } }(x)\psi ^{\sigma \Diamond _{\alpha } }(x)\Big)\\
\nonumber&\quad+\frac{1-\alpha}{{\psi ^{\Diamond _{\alpha } }(x)\psi ^{\rho\Diamond _{\alpha } }(x)}} \Big(\big(\alpha \varphi ^{\Delta\nabla }(x)+(1-\alpha )\varphi ^{\nabla \nabla }(x))\psi ^{\rho  }(x)+\psi ^{\Diamond _{\alpha } }(x)\psi ^{\nabla   }(x)\big)\psi ^{\Diamond _{\alpha } }(x)\\
\nonumber&\quad-\big(\alpha\psi ^{\Delta \nabla }(x)+(1-\alpha )\psi ^{\nabla \nabla }(x)\big)\varphi  ^{\Diamond _{\alpha } }(x)\psi(x)-\varphi ^{\nabla }(x)\psi ^{\Diamond _{\alpha } }(x)\psi ^{\rho\Diamond _{\alpha } }(x) \Big)\\
\nonumber&\qquad0\le \alpha\le1.
\end{aligned}
\end{equation}
\end{proof}

\begin{pro}\label{Proposition(3.2)}
Suppose that the function  $\psi (t)\psi ^{\sigma } (t)\psi ^{\rho } (x)\ne 0$, and $\varphi /\psi :\mathbb{T}\to R$ is diamond-alpha differentiable at $t\in\mathbb{T}$, then the following two assertions are valid.
	
(1) If
\begin{equation*}
\psi(t)\psi ^{\sigma }(t)\psi^{\rho }(x)>(<) 0
\end{equation*}
and
\begin{equation*}
\varphi ^{\Diamond _{\alpha } } (x)\psi ^{\sigma } (x)\psi ^{\rho } (x)\ge\alpha \varphi ^{\sigma }(x)\psi ^{\rho }(x)\psi ^{\Delta }(x)
+(1-\alpha )\varphi ^{\rho }(x)\psi ^{\sigma }(x)\psi ^{\nabla }(x),
\end{equation*}
then the function $\varphi /\psi$ is increasing (decreasing).
	
(2) If
\begin{equation*}
\psi(t)\psi^{\sigma }(t)\psi ^{\rho }(x)<(>) 0
\end{equation*}
and
\begin{equation*}
\varphi ^{\Diamond_{\alpha}} (x)\psi^{\sigma}(x)\psi ^{\rho}(x)\ge\alpha \varphi ^{\sigma}(x)\psi ^{\rho }(x)\psi ^{\Delta }(x)
+(1-\alpha)\varphi^{\rho}(x)\psi^{\sigma }(x)\psi^{\nabla }(x),
\end{equation*}
then the function $\varphi /\psi$ is decreasing (increasing).
	
\end{pro}

\begin{proof}
Based on Proposition \ref{Proposition(3.2)}, we know that if $(\frac{\varphi }{\psi } )^{\Diamond_{\alpha}}(x)\ge0$ which means that both
\begin{equation*}
\psi (t)\psi ^{\sigma } (t)\psi ^{\rho } (x)
\end{equation*}
and
\begin{equation*}
\varphi ^{\Diamond_{\alpha}} (x)\psi ^{\sigma } (x)\psi ^{\rho } (x)-\alpha \varphi ^{\sigma }(x)\psi ^{\rho }(x)\psi ^{\Delta }(x)-(1-\alpha )\varphi ^{\rho }(x)\psi ^{\sigma }(x)\psi ^{\nabla }(x)
\end{equation*}
are non-positive or nonnegative, then $\varphi /\psi$ is increasing.
	
If
$(\frac{\varphi }{\psi})^{\Diamond_{\alpha}}(x)\le0$ which means that both
\begin{equation*}
\psi (t)\psi ^{\sigma } (t)\psi^{\rho } (x)
\end{equation*}
is non-positive (nonnegative) and
\begin{equation*}
\varphi ^{\Diamond _{\alpha } } (x)\psi ^{\sigma } (x)\psi ^{\rho } (x)-\alpha \varphi ^{\sigma }(x)\psi ^{\rho }(x)\psi ^{\Delta }(x)-(1-\alpha )\varphi ^{\rho }(x)\psi ^{\sigma }(x)\psi ^{\nabla }(x)
\end{equation*}
is nonnegative (non-positive), then $\varphi /\psi$ is decreasing.
\end{proof}

Next, according to the definition with diamond-alpha derivatives on time scale \cite{ref6}, we give the following rule about diamond-alpha derivatives.

\begin{mr}\label{Monotonicity rule 3.1}

Suppose that the functions $\varphi,\psi$ defined on $[\beta,\gamma]_{\mathbb{T}}$ satisfy that $\varphi^{\Diamond_{\alpha}} $ and $\psi^{\Diamond_{\alpha}}$ are rd-continuous as well as $\psi^{\Diamond_{\alpha}}(x)\ne 0,\psi ^{\Diamond_{\alpha}}(x)\psi ^{\sigma \Diamond_{\alpha}} (x)\psi ^{\rho\Diamond_{\alpha}}(x)\ne0$ . If the function
\begin{equation*}
\begin{aligned}
\Big(\frac{\varphi^{\Diamond_{\alpha}}}{\psi^{\Diamond_{\alpha}}}\Big)^{\Diamond_{\alpha}}(x)
&=\frac{1}{\psi^{\Diamond_{\alpha}}(x)\psi ^{\sigma\Diamond_{\alpha}}(x)\psi ^{\rho\Diamond_{\alpha}}(x)}\Big(\varphi ^{\Diamond _{\alpha}\Diamond_{\alpha } } (x)\psi ^{\sigma\Diamond_{\alpha }} (x)\psi^{\rho\Diamond _{\alpha }}(x)\\
&\qquad-\alpha \varphi ^{\sigma\Diamond_{\alpha}}(x)\psi ^{\rho\Diamond _{\alpha}}(x)\big(\alpha\psi^{\Delta\Delta }(x)+(1-\alpha)\psi ^{\nabla\Delta}(x)\big)\\
&\qquad-(1-\alpha )\varphi ^{\rho\Diamond _{\alpha } }(x)\psi ^{\sigma\Diamond _{\alpha } }(x)\big(\psi ^{\Delta\nabla }(x)+(1-\alpha)\psi^{\nabla\nabla }(x)\big)\Big)
\end{aligned}
\end{equation*}
is nonnegative (non-positive), then the functions $\varphi^{\Diamond_{\alpha}}(x) / \psi^{\Diamond_{\alpha}}(x)$
is increasing (decreasing) on $[\beta,\gamma]_\mathbb{T}$.
\end{mr}

\end{document}